\documentclass[11pt,reqno]{amsart}

\usepackage{amssymb}
\usepackage{amsmath}
\usepackage{amscd}
\begin{document}
\setlength{\oddsidemargin}{0.5cm} \setlength{\evensidemargin}{0.5cm}

\theoremstyle{plain}
\newtheorem{theorem}{Theorem}[section]
\newtheorem{proposition}[theorem]{Proposition}
\newtheorem{lemma}[theorem]{Lemma}
\newtheorem{corollary}[theorem]{Corollary}
\newtheorem{conj}[theorem]{Conjecture}

\theoremstyle{definition}
\newtheorem{definition}[theorem]{Definition}
\newtheorem{exam}[theorem]{Example}
\newtheorem{remark}[theorem]{Remark}

\numberwithin{equation}{section}

\title[The uniqueness in the de Rham-Wu decomposition]
{The uniqueness in the de Rham-Wu decomposition}

\author{Zhiqi Chen}
\address{School of Mathematical Sciences and LPMC, Nankai University, Tianjin 300071, P.R. China} \email{chenzhiqi@nankai.edu.cn}



\begin{abstract}
In this paper, we study the uniqueness in the de Rham-Wu decomposition for pseudo-Riemannian manifolds.
\end{abstract}

\maketitle
\section{Introduction}
A {\it pseudo-Riemannian manifold}
$M$ is a smooth manifold $M$ with a nondegenerate inner product
$\langle\cdot,\cdot\rangle$ on the fibers of its tangent bundle $TM$. Let the
expression $(n_+, n_-)$, where $n_++n_-=\dim M$, denote the
signature of $\langle\cdot,\cdot\rangle$. The manifold $M$ is Riemannian if $\langle\cdot,\cdot\rangle$ has signature $(\dim M,0)$, i.e. is
positive definite. Let $H$ be the connected component of the holonomy group of a pseudo-Riemannian manifold $M$ at
some point $o$.
\medskip

There has been a lot of progress in the study on pseudo-Riemannian manifolds, such as semi-Riemannian manifolds with a doubly warped structure (\cite{GO}), the
classification of indecomposable holonomy groups of Lorentzian
manifolds (\cite{BBI,Ga1,Ga2,GT1,Le1}), pseudo-Riemannian manifolds of
index 2 (\cite{Ga3,Ik1}) and pseudo-Riemannian manifolds with neutral
signature (\cite{BBI1}) and so on. The article \cite{GT2} gives a good survey in this field.
\medskip

It is well known that one of the fundamental results in differential geometry is de Rham's decomposition theorem for Riemannian manifolds according to their holonomy representation (\cite{Be,Br1,CMS,de,EH,Pa}). The generalization to the pseudo-Riemannian case is due to Wu (\cite{Wu2,Wu1}), which is called de Rham-Wu decomposition.
But in the de Rham-Wu decomposition, the uniqueness is proved only when the maximal trivial subspace of the holonomy
group is nondegenerate. Here the maximal trivial subspace of $H$ in the tangent space $M_o$ of $M$ at $o$ is defined by
$$(M_o)^H=\{v\in M_o\mid hv=v \text{ for any } h\in H\}.$$ That nondegeneracy condition applies only to some
pseudo-Riemannian manifolds. In this paper we give a more general condition which implies uniqueness for the de Rham-Wu decomposition.
\medskip

We call $M_o$ {\it indecomposable} if it is not the
orthogonal direct sum of two proper $H$-invariant nondegenerate
subspaces, otherwise {\it decomposable}; $M_o$ {\it not reducible} if it is not the direct sum of two proper $H$-invariant subspaces,
otherwise {\it reducible}. In the Riemannian case, the term ``indecomposable" is equivalent with
the term ``not reducible". But in the pseudo-Riemannian case, there exists a pseudo-Riemannian manifold $M$ such that $M_o$ is
indecomposable but reducible. Assume that
$M_o=M_o^1\oplus\cdots\oplus M_o^p$ is an orthogonal decomposition
of $M_o$ into indecomposable $H$-invariant subspaces. Let
$\mathrm{Card}(M,o)$ denote the number of the subspaces $M_o^i$ in
the above decomposition satisfying the following conditions:
\begin{enumerate}
  \item $M_o^i$ is reducible, and
  \item $\{v\in M_o^i \mid hv=v \text{ for any } h\in H\}\not=0$.
\end{enumerate}
\medskip

An example in \cite{Wu1} shows that the decomposition of the
holonomy group into indecomposable normal subgroups is not necessary
unique up to the order if $\mathrm{Card}(M,o)=2$ for some
decomposition of $M_o$. In this paper, we prove that if $\mathrm{Card}(M,o)\leq 1$, then
\begin{enumerate}
\item $\mathrm{Card}(M,o)$ can't vary with the decomposition of $M_o$;
\item the decomposition of the tangent space $M_o$ of $M$ at $o$ into a direct sum of orthogonal subspaces
    $M_o=M_o^0\oplus M_o^1\oplus\cdots\oplus M_o^p$
  is unique up to a certain map, where $H$ acts trivially on $M_o^0$ and indecomposably on $M_o^1,\cdots,M_o^p$;
\item and the decomposition of $H$ into indecomposable normal subgroups is unique up to the order.
\end{enumerate}
Here we need to point out:
\begin{enumerate}
  \item The case given in the de Rham-Wu decomposition, i.e. the maximal trivial subspace of the holonomy
group is nondegenerate, is a special case of $\mathrm{Card}(M,o)=0$.
  \item When $\mathrm{Card}(M,o)\leq 1$, the decomposition of $M_o$ is not necessary to be unique up to the order, but the decomposition of $H$ is unique up to the order.
  \item
If $M_o$ is indecomposable but reducible, then there exists a
decomposition of $M_o$ into the direct sum of two isotropic
$H$-invariant subspaces, in particular, $M_o$ must have a neutral
signature (see Theorem~32 of \cite{BK}). What about
the inverse problem? Anyway, we have a necessary condition for $\mathrm{Card}(M,o)\leq 1$ by checking the signature of $M_o^i$. In particular, if $M$ is a Lorentzian manifold, then $\mathrm{Card}(M,o)\leq 1$. Thus the decomposition of the holonomy group of $M$ into
indecomposable normal subgroups is unique up to the order.
\end{enumerate}

\section{Preliminaries}\label{RW}

Let $S(M_o,H)$ denote the vector subspace spanned
by $v-hv$ for any $v\in M_o$ and $h\in H$.

\begin{lemma}\label{lemma1}
Let notation be as above. Then $(M_o)^H=(S(M_o,H))^\bot$.
\end{lemma}
\begin{proof}
In fact, $v\in (S(M_o,H))^\bot$, if and only if $\langle v,
w-hw\rangle=0$ for any $w\in M_o$ and $h\in H$, if and only if
$\langle v, w\rangle-\langle v,hw\rangle=0$ for any $w\in M_o$ and
$h\in H$, if and only if $\langle v, w\rangle-\langle h^{-1}v,
w\rangle=0$ for any $w\in M_o$ and $h\in H$, if and only if $\langle
v-hv, w\rangle=0$ for any $w\in M_o$ and $h\in H$, if and only if
$v-hv=0$ for any $h\in H$, i.e., $v\in (M_o)^H$.
\end{proof}

The following is the de Rham-Wu decomposition theorem.

\begin{theorem}[\cite{Wu2,Wu1}, de Rham-Wu decomposition theorem]
Let $(M,\langle\cdot,\cdot\rangle)$ be a pseudo-Riemannian manifold and let $H$ be the connected component of its holonomy group at a point $o\in M$. Then
\begin{enumerate}
   \item The tangent space $M_o$ of $M$ at $o$ decomposes into a direct sum of orthogonal subspaces
    \begin{equation}\label{1}M_o=M_o^0\oplus M_o^1\oplus\cdots\oplus M_o^p\end{equation}
such that $H$ acts trivially on $M_o^0$ and indecomposably on $M_o^1,\cdots,M_o^p$. Moreover, $H$ is a direct product
      \begin{equation}\label{2}H=H^1\times H^2\times\cdots\times H^p\end{equation}
      of normal subgroups, where each $H^i$ acts trivial on $M_o^j$, for $i\not=j$, and indecomposably on $ M_o^i$.
   \item If the maximal trivial subspace of $H$ in $M_o$ is non-degenerate with respect to the metric $\langle\cdot,\cdot\rangle$, then the decompositions in (\ref{1}) and (\ref{2}) are unique up to order.
   \item $(M,\langle\cdot,\cdot\rangle)$ is locally isometric to a pseudo-Riemannian product manifold
     \begin{equation}\label{3} (M,\langle\cdot,\cdot\rangle)\simeq (M^0,\langle\cdot,\cdot\rangle^0)\times(M^1,\langle\cdot,\cdot\rangle^1)\times\cdots\times (M^p,\langle\cdot,\cdot\rangle^p)\end{equation}
where $(M^0,\langle\cdot,\cdot\rangle^0)$ is flat or empty, and $H$ acts indecomposably on the tangent spaces of $M^i$, which
are equal to the $M_o^i$'s in (\ref{1}).
   \item If, in addition, $M$ is simply connected and $(M,\langle\cdot,\cdot\rangle)$ geodesically complete, then the decomposition
(\ref{3}) is global and the $H^i$'s in the direct decomposition (\ref{2}) are the holonomy groups of the $(M^i,\langle\cdot,\cdot\rangle)$'s.
\end{enumerate}
\end{theorem}

Obviously, if $H$ is the connected component of the holonomy group of a pseudo-Riemannian
manifold $M$ at the point $o$, then $M_o$ admits an orthogonal
decomposition into $H$-invariant subspaces
$$
M_o=M_o^0\oplus M_o^1\oplus\cdots\oplus M_o^{p_1}\oplus
M_o^{p_1+1}\oplus\cdots\oplus M_o^{p_1+p_2}, $$ where $M_o^i: 0\leq
i\leq p_1+p_2$ satisfy the following conditions:
\begin{enumerate}
  \item $M_o^0$ is a maximal nondegenerate subspace in $(M_o)^H$,
  \item $M_o^i$ is indecomposable for any $1\leq i\leq p_1+p_2$,
  \item $(M_o^i)^H=0$ for any $1\leq i\leq p_1$, and $(M_o^{p_1+i})^H\not=0$ is isotropic for any
$1\leq i\leq p_2$.
\end{enumerate}
Assume that $V$ is a nondegenerate, indecomposable and $H$-invariant
subspace with $V^H=0$. By the proof given in Appendix 1 of \cite{Wu1}, $$V=M_o^i \text{ for some }1\leq
i\leq p_1.$$ Furthermore, any orthogonal decomposition of $M_o$ into
$H$-invariant subspaces is
$$
M_o=N_o^0\oplus M_o^1\oplus\cdots\oplus M_o^{p_1}\oplus
N_o^{p_1+1}\oplus\cdots\oplus N_o^{p_1+q_2}, $$ where $N_o^0$ is a
maximal nondegenerate subspace in $(M_o)^H$, $N_o^{p_1+i}$ is
indecomposable, and $(N_o^{p_1+i})^H\not=0$ is isotropic for any
$1\leq i\leq q_2$.

\section{$\mathrm{Card}(M,o)=0$}\label{deRham}
For a pseudo-Riemannian manifold $M$, $p_2=0$ for every
decomposition of $M_o$ if $(M_o)^H$ is nondegenerate, and then
$\mathrm{Card}(M,o)=0$. The following is to discuss the case for
$p_2\geq 1$ and $\mathrm{Card}(M,o)=0$.

\subsection{De Rham decomposition of pseudo-Riemannian manifolds} This subsection is to prove the reformulated
de Rham decomposition:

\begin{theorem}\label{theorem2}
Let $M$ be a pseudo-Riemannian manifold and $H$ the connected component of its holonomy group
at the point $o$. Let
\begin{eqnarray}
M_o&=&M_o^0\oplus M_o^1\oplus\cdots\oplus M_o^{p_1}\oplus
M_o^{p_1+1}\oplus\cdots\oplus M_o^{p_1+p_2}, \label{dec1}\\
&=&N_o^0\oplus M_o^1\oplus\cdots\oplus M_o^{p_1}\oplus
N_o^{p_1+1}\oplus\cdots\oplus N_o^{p_1+q_2}
\label{dec2}\end{eqnarray} be orthogonal decompositions of $M_o$
into $H$-invariant subspaces, where
\begin{enumerate}
  \item $M_o^0$ and $N_o^0$ are maximal nondegenerate subspaces in $(M_o)^H$,
  \item $M_o^i: 1\leq i\leq p_1+p_2$ and $N_o^i: 1\leq i\leq p_1+q_2$ are indecomposable,
  \item $(M_o^i)^H=0$ for any $1\leq i\leq p_1$, $(M_o^{p_1+i})^H\not=0$ is isotropic for any $1\leq i\leq
  p_2$, and $(N_o^{p_1+j})^H\not=0$ is isotropic for any $1\leq j\leq q_2$.
\end{enumerate}
If $\mathrm{Card}(M,o)=0$ for the decomposition~(\ref{dec1}), then
we have:
\begin{enumerate}
 \item $\mathrm{Card}(M,o)=0$ for the decomposition~(\ref{dec2}),
 \item $p_2=q_2;$ $\dim M_o^{p_1+i}=\dim N_o^{p_1+i}$ and
 $S(M_o^{p_1+i},H)=S(N_o^{p_1+i},H)$ for any $1\leq i\leq p_2$ by changing the subscripts if
 necessary,
 \item For $p_1+1\leq i\leq p_1+p_2$, there exists $\pi_i$ from $M_o^i$ to $N_o^i$
 such that $\pi_i$ is $1-1$ and
$\langle\pi_i(x),\pi_i(x)\rangle=\langle x,x\rangle$ for any $x\in
M_o^i$. So
$\pi=(\pi_0,id,\cdots,id,\pi_{p_1+1},\cdots,\pi_{p_1+p_2})$ keeps the metric invariant. Here $\pi_0$ is the projection from $M_o^0$
onto $N_o^0$. That is, the decomposition is unique up to a map which keeps the metric invariant.
\end{enumerate}
\end{theorem}

It is enough to prove Theorem~\ref{theorem2} for $p_1=0$. In the
beginning, we will prove that Theorem~\ref{theorem2} holds for
$p_1=0$ and $M_o^0=0$. That is, $(M_o)^H$ is isotropic. Then $M_o$ admits an orthogonal
decomposition into $H$-invariant subspaces
\begin{equation}\label{dec3}
M_o=M_o^1\oplus\cdots\oplus M_o^p, \end{equation} where $M_o^i$ is
indecomposable and $(M_o^{i})^H\not=0$ is isotropic for any $1\leq
i\leq p$. Let
$H=H^1\times\cdots\times H^p$ be the decomposition of $H$ associated
with the decomposition~(\ref{dec3}). Let
$M_o=N_o^1\oplus\cdots\oplus N_o^q$ be another decomposition and
$H=G^1\times\cdots\times G^q$ be the decomposition of $H$
respectively.

\begin{lemma}\label{lemma3}
$S(M_o^1,H^1)\not=0$.
\end{lemma}
\begin{proof}
If not, $M_o^1=(M_o^1)^{H^1}$ by Lemma~\ref{lemma1}. Then
$M_o^1\subset (M_o)^H$, which contradicts the fact that the maximal
trivial subspace of $H$ is isotropic.
\end{proof}

Then there exists $x\in M_o^1$ and $h\in H^1$ such that
$x-hx\not=0$. Let $x=x_1+\cdots+x_q$ and $h=g_1g_2\cdots g_q$ be the
expression of $x$ and $H$ associated with
$M_o=N_o^1\oplus\cdots\oplus N_o^q$ and $H=G^1\times\cdots\times
G^q$ respectively. Then
$$0\not=x-hx=\sum_{i=1}^q(x_i-hx_i)=\sum_{i=1}^q(x_i-g_ix_i).$$ Without loss of generality,
assume that $x_1-g_1x_1\not=0$. That is, $$M_o^1\cap N_o^1\not=0$$
since $x_1-g_1x_1=x-g_1x\in M_o^1\cap N_o^1$.
\medskip

Furthermore by the assumption that $\mathrm{Card}(M,o)=0$ for the
decomposition~(\ref{dec3}), i.e., $M_o^i$ is not reducible for any
$1\leq i\leq p$, we have
\begin{lemma}\label{lemma4}
$M_o^1\cap (N_o^2\oplus\cdots\oplus N_o^q)=0$.\end{lemma}
\begin{proof}
If $M_o^1=(M_o^1\cap N_o^1)\oplus (M_o^1\cap
(N_o^2\oplus\cdots\oplus N_o^q))$, then $M_o^1\cap
(N_o^2\oplus\cdots\oplus N_o^q)=0$ by the assumption and $M_o^1\cap
N_o^1\not=0$. Or there exists $x\in M_o^1$ such that $x\not\in
(M_o^1\cap N_o^1)\oplus (M_o^1\cap (N_o^2\oplus\cdots\oplus
N_o^q))$. Let $$x=x_1+x_2$$ be the expression of $x$ associated with
the decomposition $M_o=N_o^1\oplus(N_o^2\cdots\oplus N_o^q)$. Let
$$x_1=x_1^1+x_1^2 \text{ and } x_2=x_2^1+x_2^2$$ be the expression
of $x_i, i=1,2$ associated with the decomposition
$M_o=M_o^1\oplus(M_o^2\oplus\cdots\oplus M_o^p)$. So
$x=x_1^1+x_1^2+x_2^1+x_2^2$. Thus $$x=x_1^1+x_2^1 \text{ and }
x_1^2+x_2^2=0.$$ Clearly $x_1^1-hx_1^1\in M_o^1$ and
$x_2^1-hx_2^1\in M_o^1$ for any $h\in H$. For any $h\in H$, let
$h=h_1h_2$ be the expression of $h$ associated with the
decomposition $H=H^1\times(H^2\times\cdots\times H^p)$. Then
\begin{eqnarray*}
&&x_1^1-hx_1^1=x_1^1-h_1x_1^1=x_1-x_1^2-h_1(x_1-x_1^2)=x_1-h_1x_1\in N_o^1,\\
&&x_2^1-hx_2^1=x_2^1-h_1x_2^1=x_2-x_2^2-h_1(x_2-x_2^2)=x_2-h_1x_2\in
\oplus_{i=2}^q N_o^i.
\end{eqnarray*}
It follows that for any $h\in H$,
$$hx_1^1\in x_1^1+M_o^1\cap N_o^1 \text{ and } hx_2^1\in x_2^1+M_o^1\cap (N_o^2\oplus\cdots\oplus N_o^q).$$
That is, ${\mathfrak b}_1={\mathbb R}x_1^1+M_o^1\cap N_o^1$ and
${\mathfrak b}_2={\mathbb R}x_2^1+M_o^1\cap (N_o^2\oplus\cdots\oplus
N_o^q)$ are $H$-invariant. Since $x\not\in (M_o^1\cap N_o^1)\oplus
(M_o^1\cap (N_o^2\oplus\cdots\oplus N_o^q))$ and $x=x_1^1+x_2^1$, we
have
$$(M_o^1\cap N_o^1)\oplus (M_o^1\cap (N_o^2\oplus\cdots\oplus N_o^q))\subsetneqq{\mathfrak b}_1\oplus {\mathfrak b}_2$$
If $M_o^1={\mathfrak b}_1\oplus {\mathfrak b}_2$, ${\mathfrak
b}_2=0$ by the assumption. In particular, $M_o^1\cap
(N_o^2\oplus\cdots\oplus N_o^q)=0.$ If $M_o^1\not={\mathfrak
b}_1\oplus {\mathfrak b}_2$, since $\dim M_o^1\leq \infty$,
repeating the above discussion, there exist $H$-invariant subspaces
${\mathfrak b}_1^k$ and ${\mathfrak b}_2^k$ satisfying
$$M_o^1={\mathfrak b}_1^k\oplus{\mathfrak b}_2^k.$$ By the
assumption, ${\mathfrak b}_2^k=0$. In particular, $M_o^1\cap
(N_o^2\oplus\cdots\oplus N_o^q)=0$.
\end{proof}

\begin{lemma}\label{lemma2}
The projection $\pi_1$ from $M_o^1$ to $N_o^1$ is $1-1$ and
$\langle \pi_1(x),\pi_1(x)\rangle=\langle x,x\rangle$ for any $x\in M_o^1$.
\end{lemma}
\begin{proof}
Let $\pi_1: M_o^1\rightarrow N_o^1$ be the projection from $M_o^1$
to $N_o^1$. Since $\ker\pi\subset M_o^1\cap (N_o^2\oplus\cdots\oplus
N_o^q)=0$, we have that $\pi$ is injective. For any $x\in M_o^1$,
let $x=x_1+x_2$ be the expression of $x$ associated with the
decomposition $M_o=N_o^1\oplus (N_o^2\oplus\cdots\oplus N_o^q)$.
Then $g_1x_2-x_2=0$ for any $g_1\in G^1$. Moreover for any $g_2\in
G^2\times\cdots\times G^q$,
$$x_2-g_2x_2=x-g_2x\in M_o^1\cap (N_o^2\oplus\cdots\oplus
N_o^q)=0,$$ i.e., $x_2\in (M_o)^H$. It follows that $\pi_1(M_o^1)$
is an $H$-invariant subspace of $N_o^1$. By the fact $x_2\in
(M_o)^H$, we also have $$\langle x,x\rangle=\langle
x_1+x_2,x_1+x_2\rangle=\langle x_1,x_1\rangle+\langle
x_2,x_2\rangle=\langle \pi_1(x),\pi_1(x)\rangle.$$ It follows that
$\pi_1(M_o^1)$ is an $H$-invariant nondegenerate subspace of
$N_o^1$. Then $\pi_1(M_o^1)=N_o^1$ since $N_o^1$ is indecomposable.
Namely $\pi_1$ is $1-1$.
\end{proof}

By Lemma~\ref{lemma2}, $N_o^1$ is not reducible since $M_o^1$ is
not reducible, and for any $x\in M_o^1$, $x=x_1+x_2$ where $x_1\in
N_o^1$ and $x_2\in (M_o)^H$. It follows that $S(M_o^1,H))\subset
S(N_o^1,H).$ Similarly $S(N_o^1,G^1)\subset S(M_o^1,H^1)$ since
$N_o^1$ is irreducible. That is,
$$S(M_o^1,H)=S(M_o^1,H^1)=S(N_o^1,G^1)=S(N_o^1,H).$$ Repeating the
above discussion for $j=2,3,\cdots,p$, we have proved
Theorem~\ref{theorem2} for $p_1=0$ and $M_o^0=0$.
\medskip

The following is to prove Theorem~\ref{theorem2} for $p_1=0$. That
is, $M_o$ admits an orthogonal decomposition into $H$-invariant
subspaces
\begin{equation}\label{dec5}
M_o=M_o^0\oplus M_o^1\oplus\cdots\oplus M_o^p, \end{equation} where
$M_o^0$ is a maximal nondegenerate subspace in $(M_o)^H$, $M_o^i$ is
indecomposable and $(M_o^{i})^H\not=0$ is isotropic for any $1\leq
i\leq p$. Let
$H=H^1\times\cdots\times H^p$ be the decomposition of $H$ associated
with the decomposition~(\ref{dec5}). Let $M_o=N_o^0\oplus
N_o^1\oplus\cdots\oplus N_o^q$ be another decomposition and
$H=G^1\times\cdots\times G^q$ be the decomposition of $H$
respectively. Similar to the proof of Lemma~\ref{lemma3}, we have
$$S(M_o^{1},H)\not=0.$$ Then there exists $0 \leq k\leq q$
such that $M_o^{1}\cap N_o^{k}\not=0$. Obviously $k\not=0$. So
$1\leq k\leq q$. Without loss of generality, assume that $k=1.$
Similar to the proof of Lemma~\ref{lemma4}, we have
$$
M_o^{1}\cap(N_o^0\oplus{N_o^{2}}\oplus \cdots\oplus N_o^{q})=0.
$$
\begin{lemma}\label{lemma5}
The projection $\pi_1^\prime$ from $M_o^1$ to $N_o^1$ is $1-1$.
\end{lemma}
\begin{proof}
Similar to the proof of Lemma~\ref{lemma2}, the projection
$\pi_{1}^\prime$ from $M_o^{1}$ to $N_o^{1}$ is injective, and for
any $x\in M_o^1$, $x=x_1+x_2$ where $x_1\in N_o^1$ and $x_2\in
(M_o)^H$. It follows that $\pi_{1}^\prime(M_o^1)$ is an
$H$-invariant subspace of $N_o^1$ and $$S(M_o^1,H)\subset
S(N_o^1,H).$$ Since $(M_o^1)^H$ is isotropic, by Lemma~\ref{lemma1},
we have $(M_o^1)^H\subset ((M_o^1)^H)^\perp=S(M_o^1,H).$ Let
$\{x_1,\cdots,x_r,x_{r+1},\cdots,x_n,x_{n+1}\cdots,x_{n+r}\}$ be a
basis of $M_o^1$, where $(M_o^1)^H$ is spanned by
$\{x_1,\cdots,x_r\}$, $S(M_o^1,H)$ is spanned by
$\{x_1,\cdots,x_n\}$, and the matrix of the metric associated with
the basis is
$$\left(\begin{array}{lll} 0 & 0 & I_r
\\ 0 & A_{n-r} & 0 \\ I_r & 0 & 0 \end{array}\right),$$ where $I_r$
is the identity matrix of $r\times r$ and $A_{n-r}$ is a diagonal
matrix with the element $\epsilon_i$, i.e., the sign. Then
$\{x_1,\cdots,x_r,\cdots,x_n,\pi(x_{n+1}),\cdots,\pi(x_{n+r})\}$ is
a basis of $\pi_1^\prime(M_o^1)$, and the matrix of the metric
associated with the basis is $$\left(\begin{array}{lll} 0 & 0 & I_r
\\ 0 & A_{n-r} & 0
\\ I_r & 0 & B \end{array}\right),$$ which is nondegenerate. That
is, $\pi_1^\prime(M_o^1)$ is $H$-invariant and nondegenerate. Since
$N_o^1$ is indecomposable, we have $\pi_1^\prime(M_o^1)=N_o^1$.
Namely, $\pi_1^\prime$ is $1-1$.
\end{proof}

By Lemma~\ref{lemma5}, $N_o^1$ is not reducible since $M_o^1$ is
not reducible, $\dim M_o^1=\dim N_o^1$, $(M_o^1)^H=(N_o^1)^H$, and
$S(M_o^1,H)=S(N_o^1,H)$. By the proof of Lemma~\ref{lemma5}, we have
$\pi_1^\prime(x)=x$ for any $x\in S(M_o^1,H)$ and
$$\langle\pi_1^\prime(x_i),\pi_1^\prime(x_j)\rangle=\langle x_i,x_j\rangle \text{
for any }1\leq i\leq n+r, 1\leq j\leq n.$$ Let
$x_s=x_{s_0}+x_{s_1}+x_{s_2}$ be the expression of $x_s$ associated
with the decomposition $M_o=N_o^0\oplus N_o^1\oplus
(N_o^2\oplus\cdots\oplus N_o^{q})$ for any $n+1\leq s\leq n+r$. For
any $n+1\leq , t\leq n+r$,
$$0=\langle x_s,x_t\rangle=\langle x_{s_0},x_{t_0}\rangle+\langle
x_{s_1},x_{t_1}\rangle$$ since $x_{s_2}\in
(M_o)^H\cap(N_o^2\oplus\cdots\oplus N_o^{q})$. Let
$$x_{s_1}^\prime=x_{s_1}+\frac{1}{2}\langle x_{s_0},x_{s_0}\rangle
x_{s-n}+\sum^{n+r}_{l=s+1}\langle x_{l_0},x_{l_0}\rangle x_{l-n}.$$
It is easy to check
\begin{eqnarray*}
  &\langle x_{s_1}^\prime,x_{s_1}^\prime\rangle=\langle x_{s_1},x_{s_1}\rangle + \langle x_{s_0},x_{s_0}\rangle =0, &\quad n+1\leq s\leq n+r; \\
  & \langle x_{s_1}^\prime,x_{t_1}^\prime\rangle=\langle x_{s_1},x_{t_1}\rangle + \langle x_{s_0},x_{t_0}\rangle =0, & n+1\leq s<t\leq n+r.
\end{eqnarray*}
Define $\pi_1: M_o^1\rightarrow N_o^1$ by
$$ \pi_1(x_j)=x_j, 1\leq j\leq n; \quad \pi_1(x_j)=x_{j_1}^\prime, n+1\leq j\leq n+r.$$
Then $\pi_1$ is $1-1$ from $M_o^1$ onto $N_o^1$ and
$\langle\pi(x),\pi(x)\rangle=\langle x,x\rangle$ for any $x\in
M_o^1$.
\medskip

Clearly, the projection $\pi_0: M_o^0\rightarrow N_o^0$ is $1-1$ and
$\langle \pi_0(x),\pi_0(x)\rangle=\langle x,x\rangle$ for any $x\in
M_o^0$. Repeating the above discussion for $j=2,3,\cdots,p$, we have
proved Theorem~\ref{theorem2} for $p_1=0$, i.e.,
Theorem~\ref{theorem2} holds.
\subsection{Holonomy groups of pseudo-Riemannian manifolds}\label{holonomygp}
Let $M$ be a pseudo-Riemannian manifold and $H$ the connected component of its holonomy group
at the point $o$. Then $M_o$ admits an orthogonal decomposition into
$H$-invariant subspaces:
\begin{equation}\label{dec4}
M_o=M_o^0\oplus M_o^1\oplus\cdots\oplus M_o^p, \end{equation} where
$M_o^0$ is a maximal nondegenerate subspace of $(M_o)^H$ and $M_o^i$
is indecomposable for any $1\leq i\leq p$. Let $H=H^1\times\cdots
H^p$ be the corresponding decomposition of $H$ associated with the
decomposition~(\ref{dec4}). Here $H^i$ is a normal subgroup of $H$
for any $1\leq i\leq p$, each $H^i$ is indecomposable and $H^i$ acts
trivially on $M_o^k$ if $k\not=i$.
\medskip

Let $\mathrm{Card}(M,o)=0$ for the decomposition~(\ref{dec4}). By
Theorem~\ref{theorem2}, $\mathrm{Card}(M,o)=0$ for every
decomposition, and any other orthogonal decomposition of $M_o$ into
$H$-invariant indecomposable subspaces is, by changing the order if necessary,
$$M_o=N_o^0\oplus N_o^1\oplus\cdots\oplus N_o^p,$$
where $\dim M_o^i=\dim N_o^i$ and $S(M_o^i,H)=S(N_o^i,H)$ for any
$1\leq i\leq p$. Let $H=G^1\times\cdots\times G^p$ be the
decomposition of $H$ associated with the decomposition
$M_o=N_o^0\oplus N_o^1\oplus\cdots\oplus N_o^p$. Here $G^i$ is a
normal subgroup of $H$ for any $1\leq i\leq p$, each $G^i$ is
indecomposable and $G^i$ acts trivially on $N_o^k$ if $k\not=i$.
\medskip

For any $h\in H^i$, let $h=g_1\cdots g_p$ be the expression of $h$
associated with the decomposition $H=G^1\times\cdots\times G^p$. For
any $x\in N_o^j$, by the discussion in the previous subsection,
$$x=x_1+x_2,$$ where $x_1\in M_o^j$ and $x_2\in (M_o)^H$. It follows
that, when $i\not=j$, $$hx=hx_1+hx_2=x_1+x_2=x=g_jx,$$ which shows
that $g_j$ is the identity map when $j\not=i$, i.e., $H^i\subset
G^i.$ Similarly, $G^i\subset H^i.$ Namely $H^i=G^i.$ In a word,
\begin{theorem}\label{theorem6}
Let $M$ be a pseudo-Riemannian manifold and $H$ the connected component of its holonomy group
at the point $o$. Then $H$ is the direct product of a finite number
of its normal subgroups which are indecomposable. If
$\mathrm{Card}(M,o)=0$, then the decomposition is unique up to the
order.
\end{theorem}

\section{$\mathrm{Card}(M,o)\geq 1$}\label{Lorentz}
Example 3(b) in \cite{Wu1} shows that the decomposition of the
holonomy group is not necessary unique up to the order if
$\mathrm{Card}(M,o)=2$ for some decomposition of $M_o$. For the case
$\mathrm{Card}(M,o)=1$, we have

\begin{theorem}\label{theorem3}
Let $M$ be a pseudo-Riemannian manifold and $H$ the connected component of its holonomy group
at the point $o$. Let
\begin{eqnarray}
M_o&=&M_o^0\oplus M_o^1\oplus\cdots\oplus M_o^{p_1}\oplus
M_o^{p_1+1}\oplus\cdots\oplus M_o^{p_1+p_2}, \label{d3}\\
&=&N_o^0\oplus M_o^1\oplus\cdots\oplus M_o^{p_1}\oplus
N_o^{p_1+1}\oplus\cdots\oplus N_o^{p_1+q_2}
\label{d4}\end{eqnarray} be orthogonal decompositions of $M_o$
into $H$-invariant subspaces, where
\begin{enumerate}
  \item $M_o^0$ and $N_o^0$ are maximal nondegenerate subspaces in $(M_o)^H$,
  \item $M_o^i: 1\leq i\leq p_1+p_2$ and $N_o^i: 1\leq i\leq p_1+q_2$ are indecomposable,
  \item $(M_o^i)^H=0$ for any $1\leq i\leq p_1$, $(M_o^{p_1+i})^H\not=0$ is isotropic for any $1\leq i\leq
  p_2$, and $(N_o^{p_1+j})^H\not=0$ is isotropic for any $1\leq j\leq q_2$.
\end{enumerate}
If $\mathrm{Card}(M,o)=1$ for the decomposition~(\ref{dec1}), then
we have:
\begin{enumerate}
 \item $\mathrm{Card}(M,o)=1$ for the decomposition~(\ref{dec2}),
 \item $p_2=q_2;$ $\dim M_o^{p_1+i}=\dim N_o^{p_1+i}$ and
 $S(M_o^{p_1+i},H)=S(N_o^{p_1+i},H)$ for any $1\leq i\leq p_2$ by changing the subscripts if
 necessary,
 \item For $p_1+1\leq i\leq p_1+p_2$, there exists $\pi_i$ from $M_o^i$ to $N_o^i$
 such that $\pi_i$ is $1-1$ and
$\langle\pi_i(x),\pi_i(x)\rangle=\langle x,x\rangle$ for any $x\in
M_o^i$. So
$\pi=(\pi_0,id,\cdots,id,\pi_{p_1+1},\cdots,\pi_{p_1+p_2})$ keeps the metric invariant. Here $\pi_0$ is the projection from $M_o^0$
onto $N_o^0$.
\item The decomposition of the holonomy group $H$
into indecomposable normal subgroups is unique up to the order.
\end{enumerate}
\end{theorem}

\begin{proof}
Since $\mathrm{Card}(M,o)=1$, we know only one component in the decomposition making a contribution to $\mathrm{Card}(M,o)$. Assume that it is $M_o^{p_1+p_2}$. Considering $M_o^{p_1+i}$ for $1\leq i\leq p_2-1$ and following the proof of Theorem~\ref{theorem2}, we have $N_o^{p_1+i}$ for $1\leq i\leq p_2-1$ satisfying the theorem. In order to proof this theorem, it is enough to find $N_o^{p_1+p_2}$ corresponding to $M_o^{p_1+p_2}$. Take $x\in N_o^{p_1+p_2}\cap \oplus_{i=0}^{p_1+p_2-1}M_o^i$. Let $x=x_0+x_1+\cdots+x_{p_1+p_2-1}$ be an expression of $x$ corresponding to $\oplus_{i=0}^{p_1+p_2-1}M_o^i$. Considering $hx-x$ for any $h\in H$, we have that $x_i\in (M_o^i)^H$ for any $1\leq i \leq p_1+p_2-1$ by $S(M^i_o,H)=S(N^i_o,H)$ for any $1\leq i \leq p_1+p_2-1$. That is, $x\in (M_o)^H\cap N_o^{p_1+p_2}$. Then $x$ is isotropic, which implies $x_0=0$.
Also for any $1\leq i \leq p_1+p_2-1$, $x_i\in (M^i_o)^H\subset S(M^i_o,H)=S(N^i_o,H)$. By the expression
of $x$, $x_i=0$ for any $1\leq i \leq p_1+p_2-1$. Then $x=0$. Namely $N_o^{p_1+p_2}\cap \oplus_{i=0}^{p_1+p_2-1}M_o^i=0$. It follows that the project $\pi_{p_1+p_2}$ from $N_o^{p_1+p_2}$ to $M_o^{p_1+p_2}$ is injective. By the proof of Lemma~\ref{lemma2}, $\pi_{p_1+p_2}(N_o^{p_1+p_2})=M_o^{p_1+p_2}$. Then the theorem follows.
\end{proof}

\begin{corollary}
Let $M$ be a Lorentzian manifold and $H$ the connected component of its holonomy group at the
point $o$. Then the decomposition of $H$ into its indecomposable
normal subgroups is unique up to the order.
\end{corollary}

\section{Acknowledgments}
This work is supported by National Natural Science Foundation of
China (No.11001133) and the Fundamental Research Funds for the
Central Universities. This work was completed when I visited the
University of California, Berkeley. Thanks are due to Joseph A. Wolf
for the invitation and great help to this work. I would
like to thank the anonymous referees for a great
deal of useful suggestions and comments on the contexts and the
presentation.


\begin{thebibliography}{99}
\bibitem{Be}
\textsc{M. Berger}, Sur Les groupes d'holomomie des varietes a
connexion affine et des varietes riemanniennes, \textit{Bull. Soc.
Math. France} \textbf{83} (1955), 279--330.

\bibitem{BK}
\textsc{L. Berard-Bergery} and \textsc{T. Krantz}, Representations
admitting two pairs of supplementary invariant spaces, \textit{Differential Geom. Appl.} \textbf{29} (2011), no. 1, 7--19.

\bibitem{BBI}
\textsc{L. Berard-Bergery} and \textsc{A. Ikemakhen}, On the
holonomy of Lorentzian manifolds. In \textit{Differential geometry:
Geometry in mathematical physics and related topics} (LosAngeles,
CA, 1990), Proc. Sympos. Pure Math. 54, Amer. Math. Soc.,
Providence, RI, 1993, 27--40.

\bibitem{BBI1}
\textsc{L. Berard-Bergery} and \textsc{A. Ikemakhen}, Sur
l'holonomie des varietes pseudo-riemanniennes de signature $(n,n)$,
\textit{Bull. Soc. Math. France} \textbf{125 (1)} (1997), 93--114.


\bibitem{Br1}
\textsc{R.L. Bryant}, Recent advances in the theory of Holonomy, \textit{Asterisque} \textbf{266} (2000), 351--374.

\bibitem{CMS}
\textsc{Q. Chi}, \textsc{S.A. Merkulov} and \textsc{L.J. Schwachhofer}, On the incompleteness of Berger's list of holonomy representations, \textit{Invent. Math.} \textbf{126} (1996), 391--411.

\bibitem{de}
\textsc{G. de Rham}, Sur la reductibilite d'un espace de Riemann,
\textit{Comm. Math. Helv.} \textbf{26} (1952), 329--344.

\bibitem{EH}
\textsc{J.H. Eschenburg} and \textsc{E. Heintze}, Unique decomposition of Riemannian manifolds, \textit{Proc. Amer. Math. Soc.} \textbf{126} (1998), 3075--3078.

\bibitem{Ga1}
\textsc{A.S. Galaev}, Isometry groups of Lobachevskian spaces,
similarity transformation groups of Euclidean spaces and Lorentzian
holonomy groups. \textit{Rend. Circ. Mat. Palermo (2) Suppl.}
\textbf{79} (2006), 87--97.

\bibitem{Ga2}
\textsc{A.S. Galaev}, Metrics that realize all types of Lorentzian
holonomy algebras, 2005. \textit{Int. J. Geom. Methods Mod. Phys.}
\textbf{3(5--6)} (2006), 1025--1045.

\bibitem{Ga3}
\textsc{A.S. Galaev}, Remark on holonomy groups of pseudo-Riemannian
manifolds of signature $(2,2+n)$, arXiv: math/0406397v2 [math.DG].

\bibitem{GT1}
\textsc{A.S. Galaev} and \textsc{T. Leistner}, Holonomy groups of
Lorentzian manifolds: classification, examples, and applications. In
\textit{Recent developments in pseudo-Riemannian geometry}, ESI
Lect. Math. Phys., European Math. Soc. Publ. House, Zurich 2008,
53--96.

\bibitem{GT2}
\textsc{A.S. Galaev} and \textsc{T. Leistner}, Recent developments
in pseudo-Riemannian holonomy theory, \textit{Handbook of
pseudo-Riemannian geometry and supersymmetry}, 581--627, IRMA Lect.
Math. Theor. Phys., 16, Eur. Math. Soc., Zurich, 2010.

\bibitem{GO}
\textsc{M. Guti\'errez} and \textsc{B. Olea}, Semi-Riemannian manifolds with a doubly warped structure, \textit{Rev. Mat. Iberoam.} \textbf{28(1)} (2012), 1--24.

\bibitem{Ik1}
\textsc{A. Ikemakhen}, Sur l'holonomie des varietes
pseudo-riemanniennes de signature $(2,2+n)$, \textit{Publ. Mat.}
\textbf{43 (1)} (1999), 55--84.


\bibitem{Le1}
\textsc{T. Leistner}, On the classification of Lorentzian holonomy
groups. \textit{J. Diff. Geom.} \textbf{76 (3)} (2007), 423--484.

\bibitem{Ma}
\textsc{R. Maltz}, The de Rham pruduct decomposition, \textit{J.
Diff. Geom.} \textbf{7} (1972), 161--174.


\bibitem{Pa}
\textsc{R. Pantilie}, A simple proof to the de Rham decomposition theorem, \textit{Bull. Math. Soc. Sc. Math. Roumanie} \textbf{36(84)} (1992), 341--343.



\bibitem{Wu2}
\textsc{H. Wu}, On the de Rham decomposition theorem, \textit{Ill.
J. Math.} \textbf{8} (1964), 291--311.

\bibitem{Wu1}
\textsc{H. Wu}, Holonomy groups of indefinite metrics,
\textit{Pacific J. Math.} \textbf{20} (1967), 351--392.


\end{thebibliography}
\end{document}